\documentclass{amsart}
\usepackage{amsmath}
\usepackage{amssymb}
\usepackage{enumerate}
\usepackage{enumitem}
\usepackage{mathrsfs}
\usepackage{cases}
\usepackage{color}
\usepackage{tikz}
\usepackage{appendix}
\usepackage{makecell}
\usepackage{multirow}
\usepackage{algpseudocode}
\usepackage{algorithmicx,algorithm}
\usepackage{graphicx}
\definecolor{darkblue}{rgb}{0.0, 0.0, 0.55}
\definecolor{bordeaux}{rgb}{0.34, 0.01, 0.1}
\usepackage[colorlinks,linkcolor=bordeaux,citecolor=darkblue,urlcolor=black,hypertexnames=true]{hyperref}
\usetikzlibrary{arrows,positioning}

\newtheorem{theorem}{Theorem}[section]

\newtheorem{example}[theorem]{Example}

\numberwithin{equation}{section}

\def\R{{\mathbb{R}}}

\def\N{{\mathbb{N}}}

\def\M{{\mathbf{M}}}
\def\x{{\mathbf{x}}}

\def\y{{\mathbf{y}}}

\def\a{{\boldsymbol{\alpha}}}

\def\b{{\boldsymbol{\beta}}}
\def\g{{\boldsymbol{\gamma}}}

\def\SS{{\mathscr{S}}}

\def\supp{\hbox{\rm{supp}}}

\def\int{\hbox{\rm{int}}}

\def\NF{\hbox{\rm{NF}}}

\newif\ifcomment
\commentfalse
\commenttrue
\usepackage{todonotes}

\begin{document}

\title[TSSOS: a Julia library to exploit sparsity for large-scale POPs]{TSSOS: a Julia library to exploit sparsity for large-scale polynomial optimization}
\author[V. Magron \and J. Wang]{Victor Magron \and Jie Wang}

\date{\today}

\begin{abstract}
The Julia library {\tt TSSOS} aims at helping polynomial optimizers to solve large-scale problems with sparse input data. 
The underlying algorithmic framework is based on exploiting correlative and term sparsity to obtain a new moment-SOS hierarchy involving  potentially much smaller positive semidefinite matrices. 
{\tt TSSOS} can be applied to numerous problems ranging from power networks to eigenvalue and trace optimization of noncommutative polynomials, involving up to tens of thousands of variables and constraints.
\end{abstract}

\maketitle

\section{Introduction}
The {\tt TSSOS} library intends to address large-scale polynomial optimization problems, where the polynomials in the problem's description involve only a small number of terms compared to the dense ones. 

The library is available at \href{https://github.com/wangjie212/TSSOS}{https://github.com/wangjie212/TSSOS}.
The ultimate goal is to provide semidefinite relaxations that are computationally much cheaper than those of the standard SOS-based hierarchy \cite{Las01} or its sparse version \cite{Las06,waki} based on correlative sparsity by taking into account the so-called {\em term sparsity}.
Existing algorithmic frameworks based on correlative sparsity have been implemented in the {\tt SparsePOP} solver~\cite{waki2008algorithm} and 
many applications of interest have been successfully handled, for instance certified roundoff error bounds \cite{toms18,toms17}, optimal powerflow problems~\cite{josz2018lasserre}, volume computation \cite{tacchi2019exploiting}, dynamical systems \cite{schlosser2020sparse},  noncommutative POPs~\cite{klep2021sparse},  Lipschitz constant estimation of deep networks \cite{chen2020semialgebraic,chen2021sublevel} and sparse positive definite functions \cite{mai2020sparse}.

Throughout the paper, we consider the following formulation of the polynomial optimization problem (POP):
\begin{equation}\label{pop}
(\textrm{Q}):\quad \rho^*=\inf_{\x}\,\{\,f(\x) : \x\in\mathbf{K}\,\},
\end{equation}
where the objective function $f$ is assumed to be a polynomial in $n$ variables $\x= (x_1,\ldots,x_n)$ and the feasible region $\mathbf{K}\subseteq\R^{n}$ is assumed to be defined by a finite conjunction of $m$ polynomial inequalities and $t$ polynomial equalities, namely
\begin{equation}\label{zone}
\mathbf{K} := \{\x\in\R^{n} : g_1(\x)\ge 0, \dots, g_m(\x) \geq 0, h_1(\x)=0, \dots, h_t(\x)=0 \},
\end{equation}
for some polynomials $g_1,\dots, g_m, h_1,\ldots,h_t$ in $\x$. 
For the sake of simplicity, we assume in what follows that there are only inequalities in \eqref{zone}, i.e., $t=0$.
A nowadays well-established scheme to handle $(\textrm{Q})$ is the \emph{moment-SOS hierarchy} \cite{Las01}, where SOS is the abbreviation of \emph{sum of squares}. 
The moment-SOS hierarchy provides a sequence of semidefinite programming (SDP) relaxations, whose optimal values are non-decreasing lower bounds of the global optimal value $\rho^*$ of $(\textrm{Q})$. Under certain mild conditions, the sequence of lower bounds is guaranteed to converge to $\rho^*$ generically in finite many steps. Despite of this beautiful theoretically property, the bottleneck of the scheme is that the size of the SDP relaxations become quickly intractable as $n$ increases.
Hence in order to improve the scalability, it is crucial to fully exploit the structure of POP \eqref{pop} to reduce the size of these relaxations.
A commonly present structure in large-scale polynomial optimization is {\em sparsity}. In view of this, {\tt TSSOS} implements the sparsity-adapted moment-SOS hierarchies. Here the terminology ``sparsity" is referred to the well-known correlative sparsity, or the newly proposed term sparsity \cite{chordaltssos,tssos}, or the combined both \cite{cstssos}. The idea of exploitation of sparsity behind {\tt TSSOS} is not restricted to POPs, but can be also applied to other SOS optimization problems, e.g., the computation of joint spectral radii \cite{sparsejsr} or learning of linear dynamical systems \cite{zhou2020proper,zhou2020fairness}.
{\tt TSSOS} can be also combined with efficient first-order methods to speed up the computation of the SDP relaxations themselves  \cite{mai2021constant,mai2020exploiting}.

\section{Algorithmic background and overall description}\label{sec:background}
Let $d_j=\lceil\deg(g_j)/2\rceil,j=1,\ldots,m$ and $d_{\min}=\max\{\lceil\deg(f)/2\rceil,d_1,\ldots,d_m\}$. The moment hierarchy indexed by the integer $d\ge d_{\min}$ for POP \eqref{pop} is defined by:
\begin{equation}\label{mom}
(\textrm{Q}_d):\quad\begin{cases}\inf &L_{\y}(f)\\
\textrm{s.t.}&\M_{d}(\y)\succeq0,\\
&\M_{d-d_j}(g_j\y)\succeq0,\quad j\in[m],\\
&y_{\mathbf{0}}=1.
\end{cases}
\end{equation}
Here $\M_{d}(\y)$ is the $d$-th order moment matrix, $\M_{d-d_j}(g_j\y)$ is the $(d-d_j)$-th order localizing matrix (see \cite{Las01} for more details) and $[m]:=\{1,\ldots,m\}$ for a positive integer $m$. The index $d$ is called the {\em relaxation order} of the hierarchy. Note that if $d<d_{\min}$, then $(\textrm{Q}_d)$ is infeasible.

\subsection{Correlative sparsity}\label{sec:correlative}
By exploiting correlative sparsity, we decompose the set of variables into subsets and then construct moment (localizing) matrices for each subset. Fix a relaxation order $d\ge d_{\min}$. Let $J':=\{j\in[m]\mid d_j=d\}$. For a polynomial $h=\sum_{\a}h_{\a}\x^{\a}$ ($\x^{\a}:=x_1^{\alpha_1}\cdots x_n^{\alpha_n}$), the {\em support} of $h$ is defined by $\supp(h):=\{\a\in\N^n\mid h_{\a}\ne0\}$. For $\a=(\alpha_i)\in\N^n$, let $\supp(\a):=\{i\in[n]\mid\alpha_i\ne0\}$. The {\em correlative sparsity pattern (csp) graph} associated with POP \eqref{pop} is defined to be the graph $G^{\textrm{csp}}$ with nodes $V=[n]$ and edges $E$ satisfying $\{i,j\}\in E$ if one of followings holds:
\begin{enumerate}
    \item[(i)] there exists $\a\in\supp(f)\cup\bigcup_{j\in J'}\supp(g_j)$ such that $\{i,j\}\subseteq\supp(\a)$;
    \item[(ii)] there exists $k\in[m]\backslash J'$ such that $\{i,j\}\subseteq\bigcup_{\a\in\supp(g_k)}\supp(\a)$.
\end{enumerate}
 
Let $\overline{G}^{\textrm{csp}}$ be a chordal extension of $G^{\textrm{csp}}$ and $\{I_l\}_{l\in[p]}$ be the list of maximal cliques of $\overline{G}^{\textrm{csp}}$. The polynomials $g_j, j\in[m]\backslash J'$ can be then partitioned into groups $\{g_j\mid j\in J_l\}, l\in[p]$ with respect to variable subsets $\{\x(I_l)\}_{l\in[p]}$, where $\x(I_l):=\{x_i\mid i\in I_l\}$. Consequently, we obtain the following moment hierarchy indexed by $d$ based on correlative sparsity:
\begin{equation*}\label{cs-mom}
(\textrm{Q}^{\textrm{cs}}_d):\quad\begin{cases}\inf&L_{\y}(f)\\
\textrm{s.t.}&\M_{d}(\y, I_l)\succeq0,\quad l\in[p],\\
&\M_{d-d_j}(g_j\y,I_l)\succeq0,\quad j\in J_l,l\in[p],\\
&L_{\y}(g_j)\ge0,\quad j\in J',\\
&y_{\mathbf{0}}=1.
\end{cases}
\end{equation*}
Here $\M_{d}(\y, I_l)$ ($\M_{d-d_j}(g_j\y)$) is the moment (localizing) matrix constructed with respect to the variables $\x(I_l)$.

\subsection{Term sparsity}\label{sec:term}
By exploiting term sparsity, we are able to construct moment (localizing) matrices with a block structure in an iterative manner. For $d\in\N$, let $\N^n_d:=\{\a=(\alpha_i)\in\N^n\mid\sum_{i=1}^n\alpha_i\le d\}$. Fix a relaxation order $d\ge d_{\min}$. Set $d_0:=0$ and $g_0:=1$.
Let $\SS_0 = \supp(f)\cup\bigcup_{j=1}^m\supp(g_j)\cup(2\N^n_{d})$ be the initial support. For each step $k\ge1$, let $F_{d,j}^{(k)}$ be the graph with $V(F_{d,j}^{(k)})=\N^n_{d-d_j}$ and
\begin{equation}\label{ts-eq2}
E(F_{d,j}^{(k)})=\{\{\b,\g\}\mid(\b+\g+\supp(g_j))\cap\SS_{k-1}\ne\emptyset\} 
\end{equation}
for $j\in\{0\}\cup[m]$ and let $G_{d,j}^{(k)}$ be a chordal extension of $F_{d,j}^{(k)}$. The extended support at the $k$-th step is defined as
\begin{equation*}
    \SS_{k}:=\bigcup_{i=0}^m\{\a+\b+\g\mid\a\in\supp(g_j), \{\b,\g\}\in E(G_{d,j}^{(k-1)})\textrm{ or }\b=\g\}.
\end{equation*}

\begin{example}\label{ex1}
Consider the polynomial $f=1+x_1^2+x_1x_2+x_2^2+x_1^2x_2+x_1^2x_2^2+x_2x_3+x_3^2+x_2^2x_3+x_2x_3^2+x_2^2x_3^2$. To minimize $f$ over $\R^3$, we can take the monomial basis $\{1,x_1,x_2,x_3,x_1x_2,x_2x_3\}$.
Figure \ref{tsp} shows the graph $F^{(1)}$ (without dashed edges) and its chordal extension $G^{(1)}$ (with dashed edges) for $f$, where we omit the subscripts $d,j$ since there is no constraint.

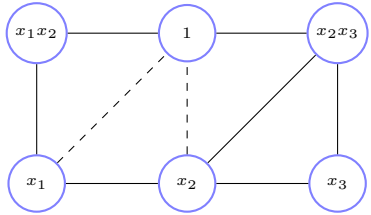
\begin{figure}[htbp]
\caption{The graph $F^{(1)}$ and its chordal extension $G^{(1)}$ for $f$}\label{tsp}
{\tiny
\begin{center}
\begin{tikzpicture}[every node/.style={circle, draw=blue!50, thick, minimum size=7.5mm}]
\node (n2) at (0,0) {$x_1$};
\node (n3) at (2,0) {$x_2$};
\node (n4) at (4,0) {$x_3$};
\node (n5) at (0,2) {$x_1x_2$};
\node (n1) at (2,2) {$1$};
\node (n6) at (4,2) {$x_2x_3$};
\draw (n2)--(n3);
\draw (n3)--(n4);
\draw (n6)--(n4);
\draw (n1)--(n6);
\draw (n1)--(n5);
\draw (n2)--(n5);
\draw (n3)--(n6);
\draw[dashed] (n2)--(n1);
\draw[dashed] (n3)--(n1);
\end{tikzpicture}
\end{center}}
\end{figure}
\end{example}

Let $\{C^{(k)}_{d,j,i}\}_{i=1}^{s_{d,j}}$ be the list of maximal cliques of $G_{d,j}^{(k)}$ for $j=0,\ldots,m$. Then the moment hierarchy based on term sparsity for POP \eqref{pop} is defined as:
\begin{equation}\label{ts-eq3}
(\textrm{Q}^{\textrm{ts}}_{d,k}):\quad
\begin{cases}
\inf &L_{\y}(f)\\
\textrm{s.t.}&[\M_d(\y)]_{C^{(k)}_{d,0,i}}\succeq0,\quad i\in[s_{d,0}],\\
&[\M_{d-d_j}(g_j\y)]_{C^{(k)}_{d,j,i}}\succeq0,\quad i\in[s_{d,j}],j\in[m],\\
&y_{\mathbf{0}}=1.
\end{cases}
\end{equation}
Here we denote by $A_C$ the submatrix of $A\in\R^{r\times r}$ with rows and columns indexed by $C\subseteq[r]$.

The above hierarchy (called the TSSOS hierarchy) is indexed by two parameters: the relaxation order $d$ and the {\em sparse order} $k$. For a fixed $d$, the sequence of optimums of $(\textrm{Q}^{\textrm{ts}}_{d,k})$ is non-decreasing and stabilizes in finitely many steps. There are two particular choices for the chordal extension $G_{d,j}^{(k)}$ of $F_{d,j}^{(k)}$: approximately smallest chordal extensions \cite{chordaltssos}\footnote{A smallest chordal extension is a chordal extension with the smallest clique number. Computing a smallest chordal extension is an NP-hard problem. Fortunately, there are efficient heuristic algorithms to produce a good approximation of smallest chordal extensions.} and the maximal chordal extension \cite{tssos}\footnote{By the maximal chordal extension, we refer to the chordal extension that completes each connected component of the graph.}. This offers a trade-off between the computational cost and the quality of obtained lower bounds. Typically, the choice of approximately smallest chordal extensions leads to positive semidefinite (PSD) blocks of smaller sizes but may provide a looser lower bound whereas the choice of the maximal chordal extension leads to PSD blocks of larger sizes but may provide a tighter lower bound. It was proved in \cite{tssos} that if the maximal chordal extension is chosen, then the sequence of optimums of $(\textrm{Q}^{\textrm{ts}}_{d,k})$ converges to the optimum of the corresponding dense relaxation as $k$ increases for a fixed relaxation order $d$.

\subsection{Correlative-term sparsity}\label{sec:cs-ts}
We can further exploit correlative sparsity and term sparsity simultaneously. Namely, first partition the set of variables into subsets $\{I_l\}_{l\in[p]}$ as done in Sec. \ref{sec:correlative} and then apply the iterative procedure for exploiting term sparsity in Sec. \ref{sec:term} to each subsystem involving variables $\x(I_l)$. Let $G_{d,l,j}^{(k)}$ be the graphs associated with each subsystem and $\{C^{(k)}_{d,l,j,i}\}_{i=1}^{s_{d,l,j}}$ be the list of maximal cliques of $G_{d,l,j}^{(k)}$ for $j\in\{0\}\cup J_l,l\in[p]$. Then the moment hierarchy based on correlative-term sparsity for POP \eqref{pop} is defined as:
\begin{equation}\label{cts}
(\textrm{Q}^{\textrm{cs-ts}}_{d,k}):\quad
\begin{cases}
\inf &L_{\y}(f)\\
\textrm{s.t.}&[\M_d(\y,I_l)]_{C^{(k)}_{d,l,0,i}}\succeq0,\quad i\in[s_{d,l,0}],l\in[p],\\
&[\M_{d-d_j}(g_j\y,I_l)]_{C^{(k)}_{d,l,j,i}}\succeq0,\quad i\in[s_{d,l,j}],j\in J_l,l\in[p],\\
&L_{\y}(g_j)\ge0,\quad j\in J',\\
&y_{\mathbf{0}}=1.
\end{cases}
\end{equation}

The above hierarchy (called the CS-TSSOS hierarchy) is also indexed by two parameters: the relaxation order $d$ and the {\em sparse order} $k$. For a fixed $d$, the sequence of optimums of $(\textrm{Q}^{\textrm{cs-ts}}_{d,k})$ is non-decreasing and stabilizes in finitely many steps. As discussed in Sec. \ref{sec:term}, the choice of chordal extensions $G_{d,l,j}^{(k)}$ offers a trade-off between the computational cost and the quality of obtained lower bounds. It was proved in \cite{cstssos} that if the maximal chordal extension is chosen, then the sequence of optimums of $(\textrm{Q}^{\textrm{cs-ts}}_{d,k})$ converges to the optimum of $(\textrm{Q}^{\textrm{cs}}_{d})$ as $k$ increases for a fixed relaxation order $d$.

\section{The implementation of {\tt TSSOS}}\label{sec:modules}
{\tt TSSOS}, which implements the sparsity-adapted moment-SOS hierarchies (i.e. $(\textrm{Q}^{\textrm{cs}}_{d})$, $(\textrm{Q}^{\textrm{ts}}_{d,k})$ and $(\textrm{Q}^{\textrm{cs-ts}}_{d,k})$), is developed as a Julia library, aiming to solve large-scale polynomial optimization problems by fully exploiting sparsity as well as other techniques. {\tt TSSOS} provides an easy way to define a POP and to solve it by SDP relaxations. The tunable parameters (e.g. $d$, $k$, the types of chordal extensions) allow the user to find the best compromise between the computational cost and the solving precision. The following script is a simple example to illustrate the usage of {\tt TSSOS}.
\vspace{3pt}

\noindent{\tt using TSSOS}\\
{\tt using DynamicPolynomials}\\
{\tt @polyvar x[1:6]}\\ 
{\tt f = x[1]\^{}4 + x[2]\^{}4 - 2x[1]\^{}2*x[2] - 2x[1] + 2x[2]*x[3] - 2x[1]\^{}2*x[3] - 2x[2]\^{}2*x[3] - 2x[2]\^{}2*x[4] - 2x[2] + 2x[1]\^{}2 + 2.5x[1]*x[2] - 2x[4] + 2x[1]*x[4] + 3x[2]\^{}2 + 2x[2]*x[5] + 2x[3]\^{}2 + 2x[3]*x[4] + 2x[4]\^{}2 + x[5]\^{}2
- 2x[5] + 2 \# define the objective function}\\
{\tt g = 1 - sum(x[1:2].\^{}2) \# define the inequality constraint}\\
{\tt h = 1 - sum(x[3:5].\^{}2) \# define the equality constraint}\\
{\tt d = 2 \# define the relaxation order}\\
{\tt numeq = 1 \# define the number of equality constraints}
\vspace{3pt}

\noindent To solve the first step of the TSSOS hierarchy with approximately smallest chordal extensions (option {\tt TS="MD"}), run
\vspace{3pt}

\noindent{\tt opt,sol,data = tssos\_first([f;g;h], x, d, numeq=numeq, TS="MD")}
\vspace{3pt}

\noindent We obtain
\vspace{3pt}

\noindent{\tt optimum = 0.20967292920706904
}
\vspace{3pt}

\noindent To solve higher steps of the TSSOS hierarchy, repeatedly run
\vspace{3pt}

\noindent{\tt opt,sol,data = tssos\_higher!(data, TS="MD")}
\vspace{3pt}

\noindent For instance, for the second step of the TSSOS hierarchy we obtain
\vspace{3pt}

\noindent{\tt optimum = 0.21230011405774876
}
\vspace{3pt}

\noindent To solve the first step of the CS-TSSOS hierarchy, run
\vspace{3pt}

\noindent{\tt opt,sol,data = cs\_tssos\_first([f;g;h], x, d, numeq=numeq, TS="MD")}
\vspace{3pt}

\noindent We obtain
\vspace{3pt}

\noindent{\tt optimum = 0.20929635879961658
}
\vspace{3pt}

\noindent To solve higher steps of the CS-TSSOS hierarchy, repeatedly run
\vspace{3pt}

\noindent{\tt opt,sol,data = cs\_tssos\_higher!(data, TS="MD")}
\vspace{3pt}

\noindent For instance, for the second step of the CS-TSSOS hierarchy we obtain
\vspace{3pt}

\noindent{\tt optimum = 0.20974835386107363
}
\vspace{3pt}

\subsection{Dependencies}\label{sec:dependencies}
{\tt TSSOS} depends on the following Julia packages:
\begin{itemize}
    \item {\tt MultivariatePolynomials} to manipulate multivariate polynomials;
    \item {\tt JuMP} \cite{jump} to model the SDP problem;
    \item {\tt LightGraphs} \cite{graph} to handle graphs;
    \item {\tt MetaGraphs} to handle weighted graphs;
    \item {\tt ChordalGraph} \cite{ChordalGraph} to generate approximately smallest chordal extensions;
    \item {\tt SemialgebraicSets} to compute Gr\"obner bases. 
\end{itemize}
Besides, {\tt TSSOS} requires an SDP solver, which can be {\tt MOSEK} \cite{mosek}, {\tt SDPT3} \cite{toh1999sdpt3}, or {\tt COSMO} \cite{cosmo}. Once one of the SDP solvers has been installed, the installation of {\tt TSSOS} is straightforward:
\vspace{3pt}

\noindent{\tt Pkg.add("https://github.com/wangjie212/TSSOS")}

\subsection{Binary variables}
{\tt TSSOS} supports binary variables. By setting ${\tt nb}=s$, one can specify that the first $s$ variables are binary variables $x_1,\dots,x_s$ which satisfy the equation $x_i^2=1$, $i \in [s]$. The specification is helpful to reduce the number of decision variables of SDP relaxations since one can identify $x^r$ with $x^{r\,(\textrm{mod }2)}$ for a binary variable $x$.

\subsection{Equality constraints}
If there are equality constraints in the description of POP \eqref{pop}, then one can reduce the number of decision variables of SDP relaxations by working in the quotient ring $\R[\x]/(h_1,\ldots,h_t)$, where $\{h_1=0,\ldots,h_t=0\}$ is the set of equality constraints. To conduct the elimination, we need to compute a Gr\"obner basis $GB$ of the ideal $(h_1,\ldots,h_t)$. Then any monomial $\x^{\a}$ can be replaced by its normal form $\NF(\x^{\a}, GB)$ with respect to the Gr\"obner basis $GB$ when constructing SDP relaxations. This reduction is conducted by default for the TSSOS hierarchy in {\tt TSSOS}.

\subsection{Adding extra first-order moment matrices}
When POP \eqref{pop} is a quadratically constrained quadratic program, the first-order moment-SOS relaxation is also known as Shor's relaxation. In this case, $(\textrm{Q}_1)$, $(\textrm{Q}^{\textrm{cs}}_1)$ and $(\textrm{Q}^{\textrm{ts}}_{1,1})$ yield the same optimum. To ensure that any higher order sparse relaxation (i.e. $(\textrm{Q}^{\textrm{cs-ts}}_{d,k})$ with $d>1$) provides a tighter lower bound compared to the one given by Shor's relaxation, we may add an extra first-order moment matrix for each variable clique in $(\textrm{Q}^{\textrm{cs-ts}}_{d,k})$:
\begin{equation}\label{cts2}
(\textrm{Q}^{\textrm{cs-ts}}_{d,k})':\quad
\begin{cases}
\inf &L_{\y}(f)\\
\textrm{s.t.}&[\M_d(\y,I_l)]_{C^{(k)}_{d,l,0,i}}\succeq0,\quad i\in[s_{d,l,0}],l\in[p],\\
&[\M_1(\y,I_l)]\succeq0,\quad l\in[p],\\
&[\M_{d-d_j}(g_j\y,I_l)]_{C^{(k)}_{d,l,j,i}}\succeq0,\quad i\in[s_{d,l,j}],j\in J_l,l\in[p],\\
&L_{\y}(g_j)\ge0,\quad j\in J',\\
&y_{\mathbf{0}}=1.
\end{cases}
\end{equation}
In {\tt TSSOS}, this is accomplished by setting ${\tt MomentOne=true}$.

\subsection{Chordal extensions}
For correlative sparsity, {\tt TSSOS} uses approximately smallest chordal extensions. For term sparsity, {\tt TSSOS} supports two types of chordal extensions: the maximal chordal extension (option ${\tt TS=``block"}$) and approximately smallest chordal extensions. {\tt TSSOS} generates approximately smallest chordal extensions via two heuristics: the Minimum Degree heuristic (option ${\tt ``TS=MD"}$) and the Minimum Fillin heuristic (option ${\tt ``TS=MF"}$). See \cite{treewidth} for a full description of these two heuristics. The Minimum Degree heuristic is slightly faster in practice, but the Minimum Fillin heuristic yields on average slightly smaller clique numbers. Hence for correlative sparsity, the Minimum Degree heuristic is recommended and for term sparsity, the Minimum Fillin heuristic is recommended.

\subsection{Merging PSD blocks}
\begin{figure}
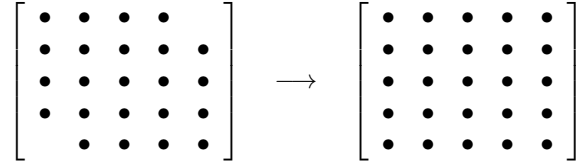

\caption{Merge two $4\times4$ blocks into a single $5\times5$ block}
\label{fig:merge}
\begin{equation*}
    \begin{bmatrix}
    \begin{array}{ccccc}
    \bullet&\bullet&\bullet&\bullet&\\
    \bullet&\bullet&\bullet&\bullet&\bullet\\
    \bullet&\bullet&\bullet&\bullet&\bullet\\
    \bullet&\bullet&\bullet&\bullet&\bullet\\
    &\bullet&\bullet&\bullet&\bullet\\
    \end{array}
    \end{bmatrix}
    \quad\longrightarrow\quad
    \begin{bmatrix}
    \begin{array}{ccccc}
    \bullet&\bullet&\bullet&\bullet&\bullet\\
    \bullet&\bullet&\bullet&\bullet&\bullet\\
    \bullet&\bullet&\bullet&\bullet&\bullet\\
    \bullet&\bullet&\bullet&\bullet&\bullet\\
    \bullet&\bullet&\bullet&\bullet&\bullet\\
    \end{array}
    \end{bmatrix}
\end{equation*}
\end{figure}

In case that two PSD blocks have a large portion of overlaps, it might be beneficial to merge these two blocks into a single block for efficiency. See Figure~\ref{fig:merge} for such an example. {\tt TSSOS} supports PSD block merging inspired by the strategy proposed in \cite{cosmo}. To activate the merging process, one just needs to set the option ${\tt Merge=True}$. The parameter ${\tt md=3}$ can be used to tune the merging strength.

\subsection{Representing polynomials in terms of supports and coefficients}
The Julia package {\tt DynamicPolynomials} provides an efficient way to define polynomials symbolically. But for large-scale polynomial optimization (say, $n>500$), it is more efficient to represent polynomials by their supports and coefficients. For instance, we can represent   $f=x_1^4+x_2^4+x_3^4+x_1x_2x_3$  in terms of its support and coefficients as follows:
\vspace{3pt}

\noindent{\tt supp = [[1; 1; 1; 1], [2; 2; 2 ;2], [3; 3; 3; 3], [1; 2; 3]] \# define the support array of f}\\
{\tt coe = [1; 1; 1; 1] \# define the coefficient vector of f}
\vspace{3pt}

The above representation of polynomials is natively supported by {\tt TSSOS}. Hence the user can define the polynomial optimization problem directly by the support data and the coefficient data to speed up the modeling process.

\subsection{Extension to noncommutative polynomial optimization}\label{sec:nctssos}
The whole framework of exploiting sparsity for (commutative) polynomial optimization can be extended to handle noncommutative polynomial optimization \cite{nctssos}, including eigenvalue and trace optimization, which leads to the submodule {\tt NCTSSOS} in {\tt TSSOS}. Table \ref{eigen_Bb} displays the numerical results for the eigenvalue minimization of the noncommutative Broyden banded function. Here ``mb'' stands for the maximal block size of the matrices involved in the SDP relaxation. 
It is evident that the sparse approach scales much better than the dense one.

\begin{table}[htbp]
\caption{The eigenvalue minimization of the noncommutative Broyden banded function: exploiting sparsity versus without exploiting sparsity. $n$: the number of variables; mb: the maixmal size of PSD blocks; opt: the optimum returned by the SDP solver; time: running time in seconds; ``-" indicates an out of memory error.}\label{eigen_Bb}
\begin{center}
\begin{tabular}{|c|c|c|c|c|c|c|}
\hline
\multirow{2}*{$n$}&\multicolumn{3}{c|}{sparse}&\multicolumn{3}{c|}{dense}\\
\cline{2-7}
&mb&opt&time&mb&opt&time\\
\hline
$20$&$15$&$0$&$0.18$&$61$&$0$&$1.39$\\
\hline
$40$&$15$&$0$&$0.72$&$121$&$0$&$66.1$\\
\hline
$60$&$15$&$0$&$1.05$&$181$&$0$&$505$\\
\hline
$80$&$15$&$0$&$1.24$&-&-&-\\
\hline
$100$&$15$&$0$&$1.41$&-&-&-\\
\hline
$200$&$15$&$0$&$3.25$&-&-&-\\
\hline
$400$&$15$&$-0.0001$&$6.70$&-&-&-\\
\hline
$600$&$15$&$-0.0002$&$13.2$&-&-&-\\
\hline
\end{tabular}
\end{center}
\end{table}


\section{Numerical experiments}
In this section, we present numerical results for the alternating current optimal power flow (AC-OPF) problem -- a famous industrial problem in power system, which can be cast as a POP involving up to tens of thousands of variables and constraints. To tackle an AC-OPF instance, we first compute a locally optimal solution with a local solver and then rely on the sparsity-adapted moment-SOS hierarchy to certify the global optimality. Suppose that the optimal value reported by the local solver is AC and the optimal value of the SDP relaxation is opt. Then the {\em optimality gap} between the locally optimal solution and the SDP relaxation is defined by
\begin{equation*}
    \textrm{gap}:=\frac{\textrm{AC}-\textrm{opt}}{\textrm{AC}}\times100\%.
\end{equation*}
When the optimality gap is less than $1\%$, we accept the locally optimal solution as globally optimal.
The test cases in this section are selected from the AC-OPF library {\em \href{https://github.com/power-grid-lib/pglib-opf}{PGLiB}} \cite{baba2019}. All numerical experiments were computed on an Intel Core i5-8265U@1.60GHz CPU with 8GB RAM memory. The SDP solver is {\tt Mosek} 9.0 with default parameters.

\begin{table}[htbp]
\caption{The results for AC-OPF instances. $n$: the number of variables; $m$: the number of constraints; AC: the local optimum; mc: the maximal size of variable cliques; opt: the optimum returned by the SDP solver; time: running time in seconds; gap: the optimality gap.}\label{opf}
\begin{center}
\small
\begin{tabular}{|c|c|c|c|c|c|c|c|}
\hline
case&$n$&$m$&AC&mc&opt&time&gap\\
\hline
14\_ieee\_api&38&147&5.9994e3&6&5.9994e3&0.54&0.00\%\\
\hline
30\_ieee&72&297&8.2085e3&8&8.2085e3&0.99&0.00\%\\
\hline
39\_epri\_sad&98&361&1.4834e5&8&1.4831e5&1.45&0.02\%\\
\hline
118\_ieee&344&1325&9.7214e4&10&9.7214e4&7.71&0.00\%\\
\hline
179\_goc\_api&416&1827&1.9320e6&10&1.9226e6&9.69&0.48\%\\
\hline
300\_ieee&738&2983&5.6522e5&14&5.6522e5&25.2&0.00\%\\ 
\hline
793\_goc&1780&7019&2.6020e5&18&2.5932e5&66.1&0.34\%\\
\hline
1354\_pegase\_sad&3228&13901&1.2588e6&26&1.2582e6&387&0.05\%\\
\hline
1951\_rte\_api&4634&18921&2.4108e6&26&2.4029e6&596&0.32\%\\
\hline
2000\_goc\_api&4476&23009&1.4686e6&42&1.4610e6&1094&0.51\%\\
\hline
2312\_goc&5076&21753&4.4133e5&68&4.3858e5&997&0.62\%\\
\hline
3022\_goc\_sad&6698&29283&6.0143e5&50&5.9859e5&1340&0.47\%\\
\hline
\end{tabular}
\end{center}
\end{table}

From Table \ref{opf}, we can see that for all test cases, {\tt TSSOS} successfully reduces the optimality gap to less than $1\%$, namely, certifies the global optimality. The largest instance 3022\_goc\_sad has 6698 variables and 29283 constraints.
\vspace{1em}

\paragraph{\textbf{Acknowledgements}.} 
This work was supported by the Tremplin ERC Stg Grant ANR-18-ERC2-0004-01 (T-COPS project).
The firts author was supported by the FMJH Program PGMO (EPICS project) and  EDF, Thales, Orange et Criteo.
This work has benefited from  the European Union's Horizon 2020 research and innovation program under the Marie Sklodowska-Curie Actions, grant agreement 813211 (POEMA) as well as from the AI Interdisciplinary Institute ANITI funding, through the French ``Investing for the Future PIA3'' program under the Grant agreement n$^{\circ}$ANR-19-PI3A-0004.


\end{document}